\documentclass[oneside,english]{amsart}
\usepackage[T1]{fontenc}
\usepackage[latin1]{inputenc}
\usepackage{amssymb}

\makeatletter

 \theoremstyle{plain} 
 \newtheorem{thm}{Theorem}[section]
 \numberwithin{equation}{section} 
 \numberwithin{figure}{section} 
 \theoremstyle{plain}
 \theoremstyle{plain}    
 \newtheorem{lem}[thm]{Lemma} 
 \theoremstyle{definition}
 \newtheorem{defn}[thm]{Definition}
 \theoremstyle{plain}    
 \newtheorem{prop}[thm]{Proposition} 
 \newtheorem{rem}[thm]{Remark}
 \theoremstyle{remark}
 \newtheorem*{rem*}{Remark}
 \theoremstyle{plain}    

 \newtheorem*{cor*}{Corollary}
 \newtheorem{cor}[thm]{Corollary}
 \theoremstyle{definition}
  \newtheorem{example}[thm]{Example}
\usepackage[T1]{fontenc}
\usepackage[latin1]{inputenc}

\makeatletter

\usepackage[T1]{fontenc}
\usepackage[latin1]{inputenc}

\makeatletter

\usepackage[T1]{fontenc}
\usepackage[latin1]{inputenc}

\makeatletter

\usepackage[T1]{fontenc}
\usepackage[latin1]{inputenc}

\makeatletter

\usepackage[T1]{fontenc}
\usepackage[latin1]{inputenc}


\makeatother

\makeatother

\makeatother

\usepackage{babel}
\makeatother

\def\lra{\longrightarrow}
\def\R{{\mathbb R}}

\def\oH1{\buildrel\circ\over H\kern-.04in{}^1}

\def\barB{{\overline{B}}}
\def\barD{{\overline{D}}}
\def\barv{{\overline{v}}}

\def\rharp{{\rightharpoonup}}

\def\meas{{\,meas\,}}
\def\supp{{\,supp\,}}
\def\const{{\,const\,}}
\def\loc{{\,loc\,}}

\def\bee{\begin{equation*}}
\def\eee{\end{equation*}}
\def\be{\begin{equation}}
\def\ee{\end{equation}}

\def\ep{{\epsilon}}

\begin{document}

\title{Embedding operators and boundary-value problems for
rough domains
\footnote{This work was supported by the
Center of Advanced Studies in Mathematics at Ben Gurion
University }}

\author{V Goldshtein$^{\textrm{*}}$ and A.G. Ramm$^{\textrm{**}}$}

\address{$^{\textrm{*}}$Mathematics Department, Ben Gurion University of
the Negev \\
 P.O.Box 653, Beer Sheva, 84105, Israel\\
 vladimir@bgumail.bgu.ac.il
}

\address{$^{\textrm{**}}$Mathematics Department, Kansas State University
\\
 Manhattan, KS 66506-2602, USA\\
 ramm@math.ksu.edu}

\begin{abstract}

In the first part of the paper boundary-value problems are considered under
weak assumptions on the smoothness of the domains.  
We assume nothing about smoothness of the boundary $\partial D$
of a bounded domain $D$ when the homogeneous Dirichlet boundary condition is
imposed; we assume boundedness of the embedding $i_{1}:H^{1}(D)\rightarrow
L^{2}(D)$ when the Neumann boundary condition is imposed; we assume
boundedness of the embeddings $i_{1}$ and of $i_{2}:H^{1}(D)\rightarrow
L^{2}(\partial D)$ when the Robin boundary condition is imposed, and, if, in
addition, $i_{1}$ and $i_{2}$ are compact, then the boundary-value problems
with the spectral parameter are of Fredholm type. Here $i_{1}$ is the
embedding of $H^{1}(D)$ (or $H^{1}(\widetilde{D})$) into $L^{2}(D)$
$(L^{2}(\widetilde{D}))$, and $\widetilde{D}$ is such bounded subdomain of
the exterior domain $D':=\R^{n}\setminus D$ that its boundary consists of
two components:  $\partial D$ and
$S$, where $S$ is a smooth compact manifold. The space $L^{2}(\partial
D)$ is the $L^{2}$ space on $\partial D$ with respect to Hausdorff
$(n-1)\hbox{-dimensional}$ measure on $\partial D$. These results
motivate our detailed study of the embedding operators.  In the second 
part of
the paper new classes of rough bounded domains $D$  are introduced.
The
embedding operator from $H^{1}(D)$ into $L^{2}(\partial D)$ is compact
for these classes of domains.
These classes include, in particular, the domains whose boundaries
are compact Lipschitz manifolds, but much larger sets of
domains are also included in the above classes. Several examples of 
the classes of rough domains for which the 
 embedding $i_2$ is
compact are given. Applications to scattering by rough obstacles are
mentioned.

\end{abstract}
\maketitle

\section{Introduction}

Embedding inequalities are studied in this paper for rough, that is, 
non-smooth, domains.
An essentially self-contained presentation of a method for a
study of boundary-value problems for second-order elliptic
equations in such domains is developed. 
The novel points include the usage of the
limiting absorption principle for the proof of the existence
of solutions and weaker than usual assumptions on the
smoothness of the boundary. For brevity of the presentation
we consider the boundary-value problems for Laplacian, and
the three classical boundary conditions. We study interior
and exterior boundary-value problems and obtain the
existence results and the Fredholm property under weak
assumptions on the smoothness of the boundary. The method we
use is applicable for general second-order elliptic
equations and for obstacle scattering problems. Elliptic boundary-value 
problems were studied in
numerous books and papers. We mention \cite{GT} and
\cite{LU}, where many references can be found. In \cite{Maz}
embedding theorems for a variety of  non-smooth domains have
been studied. In \cite{r409} the obstacle scattering
problems were studied for non-smooth obstacles. In \cite{r190} and
\cite{r470} the boundary-value problems and direct and
inverse obstacle scattering problems have been studied. In
\cite{GR} embedding theorems in some classes of
rough domains were studied. The aim of this paper is to
continue the studies initiated in \cite{r470}, \cite{r409} 
and
\cite{GR}.

Consider the boundary-value problems \be\label{e1.1} -\Delta u=F
\hbox{\ in\ } D,\quad F\in L^2(D), \ee
\be\label{e1.2} u=0\hbox{\ on\ }\partial D. \ee
The boundary conditions can be the Neumann one \be\label{e1.3} u_N=0
\hbox{\ on\ } \partial D, \ee
where $N$ is the outer unit normal to $\partial D$, or the Robin one: \be\label{e1.4}
u_N+h(s)u=0 \hbox{\ on\ }\partial D, \ee
where $h(s)\geq0$ is a bounded piecewise-continuous function on $\partial D$.

We are interested in similar problems in the exterior domain $D':=R^{n}\setminus D$,
and we consider the case $n=3$. The case $n>3$ can be treated similarly.
If $n=2$ some additional remarks are in order since the fundamental
solution in this case changes sign and tends to infinity as $|x-y|:=r_{xy}\rightarrow\infty$.
If $n=3$, then $g(x,y):=\frac{1}{4\pi r_{xy}}$, and if $n=2$, then
$g(x,y)=\frac{1}{2\pi}\ln\frac{1}{r_{xy}}$, $x,y\in\R^{n}$, $-\Delta g=\delta(x-y)$
in $\R^{n}$, and $\delta(x)$ is the delta-function.

Below $(\cdot,\cdot)$ denotes the inner product in $L^{2}(D):=H^{0}$,
$L_{0}^{2}(D)$ is the set of $L^{2}(D)$ functions with compact support,
$L_{0}^{2}(D')$ is the set of $L^{2}(D')$ functions which vanish
near infinity, $\oH1$ is the closure of $C_{0}^{\infty}(D)$ in the
$H^{1}:=H^{1}(D)$ norm $\Vert u\Vert_{1}:=\left(\int_{D}(|u|^{2}+|\nabla u|^{2})dx\right)^{1/2}$.
We denote $\Vert u\Vert:=\left(\int_{D}|u|^{2}dx\right)^{1/2}$.

If the boundary conditions are non-homogeneous, e.g., $u=f$ on $\partial D$,
then we assume that there exists a function $v\in H^{1}(D)\bigcap H_{\loc}^{2}(D)$,
$\Delta v\in L^{2}(D)$, such that $v=f$ on $\partial D$ and consider $w:=u-v$.
The function $w$ satisfies equation \eqref{e1.1} with $F$ replaced
by $F+\Delta v$, and $w$ satisfies \eqref{e1.2}. Similarly one
treats inhomogeneous Neumann and Robin boundary conditions. In the
case of inhomogeneous boundary conditions the smoothness assumptions
on the boundary $\partial D$ are more restrictive than in the case of the
homogeneous boundary conditions.

Let us reformulate the problems \eqref{e1.1}--\eqref{e1.4} so that
the assumptions on $\partial D$ are minimal.

In the case of the Dirichlet problem \eqref{e1.1}--\eqref{e1.2}
we use the weak formulation:

$u$ solves \eqref{e1.1}--\eqref{e1.2} iff $u\in\oH1(D):=\oH1$
and \be\label{e1.5} {[}u,\phi{]}:=(\nabla u,\nabla\phi)=(F,\phi)
\qquad \forall \phi \in\oH1. \ee

The weak formulation \eqref{e1.5} of the Dirichlet problem \eqref{e1.1}--\eqref{e1.2}
does not require any smoothness of $\partial D$.

The weak formulation of the Neumann problem \eqref{e1.1}, \eqref{e1.3}
is: \be\label{e1.6} {[}u,\phi{]}=(F,\phi) \qquad \forall \phi\in H^1.
\ee
An obvious necessary condition on $F$ for \eqref{e1.6} to hold is
\be\label{e1.7} (F,1)=0. \ee
Although the statement of the problem \eqref{e1.6} does not require
any smoothness assumption on $\partial D$, one has to assume that $\partial D$ is
smooth enough for the Poincare-type inequality to hold: \be\label{e1.8}
inf_{m\in\R^1} \|u-m\|\leq c\|\nabla u\|, \qquad c=\const>0,
\ee
see Remark 2.1 below.

The infimum in \eqref{e1.8} is attained at $m_{0}=\frac{1}{|D|}\int_{D}udx$,
$|D|:=\meas D$, and $(u-m_{0},1)=0$. If $(u,1)=0$, then \eqref{e1.8}
implies $\Vert u\Vert\leq c\Vert\nabla u\Vert$. The role of this
inequality will be clear from the proof of the existence of the solution
to \eqref{e1.6} (see Section 2).

Finally, for the Robin boundary condition the weak formulation of
the boundary-value problem \eqref{e1.1}, \eqref{e1.4} is: \be\label{e1.9}
{[}u,\phi{]}+\int_{\partial D} hu \bar{\phi} ds=(F,\phi) \quad \forall \phi\in 
H^1. \ee
For \eqref{e1.9} to make sense, one has to be able to define $u$
on $\partial D$. For this reason we assume that the embedding $i_{2}:H^{1}(D)\rightarrow L^{2}(\partial D)$
is bounded. We also assume the compactness of $i_{2}$, and this assumption
is motivated in the proof of the existence and uniqueness of the solution
to \eqref{e1.9}.

Let us formulate our results. We assume that $D\subset\R^{n}$, ($n=3$),
is a bounded domain and $F\in L^{2}(D)$ is compactly supported. This
assumption will be relaxed in Remark 3.1.

\begin{thm}
\label{T:1.1} The solution $u\in\oH1(D)$ of \eqref{e1.5} exists
and is unique.
\end{thm} 

\begin{thm}
\label{T:1.2} If $D$ is such that \eqref{e1.8} holds and $F$ satisfies
\eqref{e1.7}, then there exists a solution $u$ to \eqref{e1.6},
and $\{ u+c\}$, $c=\const$, is the set of all solutions to \eqref{e1.6}
in $H^{1}$. 
\end{thm}

\begin{thm}
\label{T:1.3} If $D$ is such that $i_{1}:H^{1}(D)\rightarrow L^{2}(D)$
is compact and $i_{2}:H^{1}(D)\rightarrow L^{2}(\partial D)$ is bounded,
$F\in L_{0}^{2}(D)$ and $h\geq0$ is a piecewise-continuous bounded
function on $\partial D$, $h\not\equiv0$, then problem \eqref{e1.9} has
a solution in $H^{1}(D)$ and this solution is unique. If $i_{1}$
and $i_{2}$ are compact, then the problem \[
[u,\phi]+\int_{\partial D}hu\bar{\phi}ds-\lambda(u,\phi)=(F,\phi),\quad\lambda=\const\in\R\]
 is of Fredholm type. 
\end{thm}

Similar results are obtained in Section 3 for the boundary-value problems
in the exterior domains (Theorem 3.1).
 
   Theorems \ref{T:1.1}-\ref{T:1.3} demonstrate existence of a connection between elliptic boundary problems and compactness of embedding operators for Sobolev spaces. This is a motivation for a detailed study of embedding operators.

In Sections 4 and 5 of the paper we prove some results about compactness
of the embedding operator from $H^{1}(D)$ to $L^{2}(\partial D)$
for rough bounded domains.
First, we prove compactness of the embedding operators for {}``elementary''
domains whose boundaries are Lipschitz manifolds or
even Lipschitz manifolds {}``almost everywhere'' (in some special
sense). This class of {}``elementary'' domains is larger than the known
classes of domains used for embedding theorems.

Using a lemma for the union of {}``elementary'' domains
we extend this result to domains of the class $Q$ which consists
of the finite
unions of the ``elementary'' domains. We show by  examples 
that the boundary of a bounded domain of class $Q$ can have countably
many connected components (see example \ref{component}). Another
example demonstrate that boundaries of such domains are not necessarily
have local presentation as graphs of Lipschitz functions 
(see example \ref{spiral}).

\section{Proofs of Theorems \ref{T:1.1}-\ref{T:1.3}}

\begin{proof} {[}Proof of Theorem \ref{T:1.1} {]} One has \be\label{e2.1}
|{[}u,\phi{]}|=|(F,\phi)| \leq \|F\| \|\phi\| \leq c\|F\| \|\varphi\|_1 \ee
where we have used the inequality \be\label{e2.2} \|\phi\|\leq c\|\phi\|_1,
\quad \phi\in\oH1, \ee which holds for any bounded domain, i.e., without any
smoothness assumptions on $D$. Note that the norm $[u,u]^{1/2}:=[u]$ is
equivalent to $H^{1}$ norm on $\oH1:c_{1}\Vert u\Vert_{1}\leq[u]\leq\Vert
u\Vert_{1}$, $c_{1}=\const>0$. Inequality \eqref{e2.1} shows that $(F,\phi)$
is a bounded linear functional in $H^{1}(D)$ so, by the Riesz theorem about
linear functionals in a Hilbert space, one has \[
[u,\phi]=[BF,\phi]\quad\forall\phi\in\oH1,\]
 where $B$ is a bounded linear operator from $L^{2}(D)$ into $\oH1$.
Thus $u=BF$ is the unique solution to \eqref{e1.5}.
\end{proof}

\begin{proof} {[}Proof of Theorem \ref{T:1.2}{]} If $(F,1)=0$, then one may
assume that $(\phi,1)=0$ because $(F,\phi)=(F,\phi-m)$ and the constant $m$
can be chosen so that $(\phi-m,1)=0$ if $D$ is bounded. If $D$ is such that
\eqref{e1.8} holds, then \be\label{e2.3} |{[}u,\phi{]}| = |(F,\phi-m)| \leq
\|F\| \inf_m \|\phi-m\| \leq c\|F\| \|\nabla\phi\|. \ee Thus
$(F,\phi)=[BF,\phi]$, where $B:L^{2}(D)\rightarrow H^{1}$ is a bounded
linear operator. Thus $u=BF$ solves \eqref{e1.6}, for any constant $m$ and
$u+m$ also solves \eqref{e1.6} because $[m,\phi]=0$. If $u$ and $v$ solve
\eqref{e1.6}, then $w:=u-v$ solves the equation
$[w,\phi]=0\quad\forall\phi\in H^{1}$. Take $\phi=w$ and get
$[w,w]=\Vert\nabla w\Vert^{2}=0$. Thus $\nabla w=0$ and $w=\const$. Theorem
\ref{T:1.2} is proved. \end{proof}

\begin{rem}
Necessary and sufficient conditions on \(D\) for \eqref{e1.8} to hold one
can find in \cite{Maz}. Inequality \eqref{e1.8} is equivalent to the
boundedness of the embedding in \(i_1:L^1_2(D)\to L^2(D)\)  
\cite[p.169]{Maz}.
By \(L^1_2\) denote the space of functions \(u\) such that \(\|\nabla u\|<\infty\)
is denoted.
\end{rem}

\begin{rem}
If one wants to study the problem
\be\label{e2.4}
  -\Delta u-\lambda u=F, \qquad u=0\hbox{\ on\ } \partial D \ee where
\(\lambda=\const\), and a similar problem with the Neumann boundary
condition \eqref{e1.3} or with the Robin condition \eqref{e1.4} to be of
Fredholm type, then one has to assume the operators \(B\) in the proofs of
Theorem \ref{T:1.1} and Theorem \ref{T:1.2} to be compact in \(\oH1\) and in
\(H^1\) respectively. Originally the operators \(B\) were acting from
\(L^2(D)\) onto \(\oH1\) and \(H^1\) (respectively in Theorem \ref{T:1.1}
and in Theorem \ref{T:1.2}). Thus, \(B\) are defined on \(\oH1\subset
L^2(D)\) and on \(H^1\subset L^2(D)\) respectively. \end{rem}

\begin{proof} {[}Proof of Theorem 1.3.{]} If the embedding
$i_{2}:H^{1}(D)\rightarrow L^{2}(\partial D)$ is bounded, then
\be\label{e2.5} |\int_{\partial D} hu\bar{\phi} ds| \leq \sup_{\partial
D}|h| \|u\|_{L^2(\partial D)} \|\phi\|_{L^2(\partial D)} \leq c\|\phi\|_1.
\ee By Riesz's theorem one gets \[ \int_{\partial
D}hu\bar{\phi}ds=(Tu,\phi)_{1}.\]
 Equation \eqref{e1.9} can be written as \be\label{e2.6} (Au+Tu-BF,
\phi)_1=0 \qquad \forall \phi\in H^1, \ee where $[u,\phi]=(Au,\phi)_{1}$,
$(F,\phi)=(BF,\phi)_{1}$. Thus \be\label{e2.7} Qu-BF:= Au+Tu-BF=0, \ee where
$A$ is a bounded linear operator in $H^{1}$, $\Vert A\Vert\leq1$ because
$[u,u]\leq(u,u)_{1}$, $||B||\leq1$ because
$|(BF,\phi)_{1}|=|(F,\phi)|\leq\Vert F\Vert\Vert\phi\Vert\leq\Vert
F\Vert_{1}\Vert\phi\Vert_{1}$, and $T$ is a bounded operator in $H^{1}$ if
the embedding operator $i_{2}:H^{1}(D)\rightarrow L^{2}(\partial D)$ is
bounded. If $i_{2}$ is compact, then $T$ is compact in $H^{1}$. The operator
$Q:=A+T$ is linear, defined on all of $H^{1}$, and bounded. The expression
\[ N^{2}(u):=(Qu,u)_{1}=[u,u]+\int_{\partial D}h|u|^{2}ds\]
 defines a norm $N(u)$ equivalent to $\Vert u\Vert_{1}$.

Let us prove this equivalence.

By \eqref{e2.5} one has $N^{2}(u)\leq c\Vert u\Vert_{1}^{2}$. Also
\[
\Vert u\Vert_{1}^{2}=[u,u]+(u,u)\leq N^{2}(u)+(u,u)\leq c_{1}N^{2}(u)\]
 because \[
\Vert u\Vert\leq cN(u),\]
 where $c=\const>0$ stands for various constants independent of $u$.

\textit{Let us prove the inequality $\Vert u\Vert\leq cN(u)$.}

Assuming that it fails, one finds a sequence $u_{n}\in
H^{1}$, $\Vert u_{n}\Vert=1$, such that $\Vert
u_{n}\Vert\geq nN(u_{n})$, so $N(u_{n})\leq\frac{1}{n}$.
Thus $\Vert\nabla u_{n}\Vert\rightarrow0$ and
$\int_{\partial D}h|u_{n}|^{2}ds\rightarrow0$. Since $\Vert
u_{n}\Vert=1$ one may assume that $u_{n}\rharp v$, where
$\rharp$ denotes weak convergence in $L^{2}(D)$. If
$u_{n}\rharp v$ and $\nabla u_{n}\rharp0$, then $\nabla
v=0$, so $v=C=\const$, and \[
0=\lim_{n\rightarrow\infty}\int_{\partial
D}h|u_{n}|^{2}ds=C^{2}\int_{\partial D}hds,\]
 so $C=0$. Since the embedding $i_{1}:H^{1}(D)\rightarrow
L^{2}(D)$ is compact and the sequence $u_{n}$ is bounded in
$H^{1}$, we may assume without loss of generality that
$u_{n}$ converges in $L^{2}(D)$ to zero. This contradicts
the assumption $\Vert u_{n}\Vert=1$.  The inequality is
proved.

Thus, the norms $N(u)$ and $\Vert u\Vert_{1}$ are equivalent, the operator
$Q$ is positive definite, selfadjoint as an operator in $H^{1},$ and
therefore $Q$ has a bounded inverse in $H^{1}$. Thus, equation \eqref{e2.7}
has a unique solution $u=Q^{-1}BF$ in $H^{1}$. The statement of Theorem
\ref{T:1.3} concerning Fredholm's type of problem \eqref{e1.9} follows from
Lemma \ref{L:1} below. Theorem \ref{T:1.3} is proved \end{proof}

\begin{rem}
As in Remark 2.2, if \(B\) is compact in \(H^1\), then the problem
\be\label{e2.8}
  [u,\phi]+\int_{\partial D} hu\bar{\phi} ds=\lambda(u,\phi)+(F,\phi),
  \quad \lambda=\const \ee
is of Fredholm type. This problem can be written as \(Au+Tu=\lambda Bu+BF\), or
\be\label{e2.9}
  u=\lambda Q^{-1} Bu+Q^{-1} BF \ee
where the operator \(Q^{-1}B\)  is compact in \(H^1\).
\end{rem}

\begin{lem}
\label{L:1} The operator $B$ is compact in $H^{1}$ if and only
if the embedding operator $i_{1}:H^{1}(D)\rightarrow L^{2}(D)$
is compact. 
\end{lem}

\begin{proof}
Suppose that the embedding 
$i_{1}:H^{1}(D)\rightarrow L^{2}(D)$ is compact.
One has $\Vert u\Vert^{2}=(Bu,u)_{1}=(u,Bu)_{1}$, so $B$ is a linear
positive, symmetric, and bounded operator in $H^{1}$. One has
$(u,\phi)=(Bu,\phi)_{1}$, so $\Vert Bu\Vert_{1}\leq\Vert u\Vert\leq\Vert
u\Vert_{1}$, so $\Vert B\Vert_{H^{1}\rightarrow H^{1}}\leq1$. A linear
positive, symmetric, bounded operator $B$, defined on all of $H^{1}$, is
selfadjoint. The operator $B^{1/2}>0$ is well defined, $B$ and $B^{1/2}$ are
simultaneously compact, and $\Vert u\Vert=\Vert B^{1/2}u\Vert_{1}$. Thus, if
$j$ is compact then the inequality $\Vert u_{n}\Vert_{1}\leq1$ implies the
existence of a convergent in $L^{2}(D)$ subsequence, denoted again $u_{n}$,
so that $B^{1/2}u_{n}$ converges in $H^{1}$. Thus, $B^{1/2}$ is compact in
$H^{1}$ and so is $B$.

Conversely, if $B$ is compact in $H^{1}$ so is $B^{1/2}$. Therefore, if
$u_{n}$ is a bounded in $H^{1}$ sequence, $\Vert u_{n}\Vert_{1}\leq1$, then
$B^{1/2}u_{n_{k}}$ is a convergent in $H^{1}$ sequence. Denote the
subsequence $u_{n_{k}}$ again $u_{n}$. Then $u_{n}$ is a convergent in
$H^{0}=L^{2}(D)$ sequence because $\Vert u_{n}\Vert=\Vert
B^{1/2}u_{n}\Vert_{1}$. Therefore $i_{1}$ is compact. Lemma \ref{L:1} is
proved. \end{proof}

\begin{rem}
We have used the assumptions \(h\geq 0\) and \(h\not\equiv 0\) in the proof
of Theorem \ref{T:1.3} .
If \(h\) changes sign on \(\partial D\) but the embeddings \(i_2:H^1(D)\to 
L^2(\partial D)\)
and \(i_1:H^1(D)\to L^2(D)\) are compact,
then problem \eqref{e1.9} is still of Fredholm's type because
\(T\) is compact in \(H^1\) if \(i_2\) is compact.
\end{rem}

\section{Exterior boundary-value problems}

Consider boundary-value problems \eqref{e1.1}, \eqref{e1.2}, \eqref{e1.3},
\eqref{e1.4} in the exterior domain $D'=\R^{3}\setminus\barD$. The closure
of $H_{0}^{1}(D')$ in the norm $\Vert
u\Vert_{1}:=\{\int_{D'}(|u|^{2}+|\nabla u|^{2})dx\}^{1/2}$ is denoted by
$H^{1}=H^{1}(D')$ and $H_{0}^{1}(D')$ is the set of functions vanishing near
infinity and with finite norm $\Vert u\Vert_{1}<\infty$. We assume that $D$
is bounded. The weak formulation of the boundary-value problems is given
similarly to \eqref{e1.5}, \eqref{e1.6} and \eqref{e1.9}. The corresponding
quadratic forms Dirichlet $t_{D}$, Neumann $t_{N}$ and Robin $t_{R}$, where
$t_{D}[u,u]=(\nabla u,\nabla u)$, $u\in\oH1(D')$; $t_{N}[u,u]=(\nabla
u,\nabla u)$, $u\in H^{1}(D')$; $t_{R}[u,u]=(\nabla u,\nabla u)+\langle
hu,u\rangle$, $u\in H^{1}(D')$, $\langle u,v\rangle:=\int_{\partial D}u\barv
ds$, $0\leq h\leq c$, are nonnegative, symmetric and closable. Here and
below, $c>0$ stands for various constants. Nonnegativity and symmetry of the
above forms are obvious.

\textit{Let us prove their closability.}

By definition, a quadratic form $t[u,u]$ bounded from below, i.e.,
$t[u,u]>-m(u,u)$, $m=\const$, and densely defined in the Hilbert space
$H=L^{2}(D')$, is closable if
$t[u_{n}-u_{m},u_{n}-u_{m}]\underset{n,m\rightarrow\infty}{\lra}0$ and
$u_{n}\underset{H}{\lra}0$ imply
$t[u_{n},u_{n}]\underset{n\rightarrow\infty}{\lra}0$. The closure of the
domain $D[t]$ of the closable quadratic form in the norm $[u]:=\{
t[u,u]+(m+1)(u,u)\}^{1/2}$ is a Hilbert space $H_{t}\subset H$ densely
embedded in $H$. The quadratic form $t[u,u]$ is defined on $H_{t}$ and this
form with the domain of definition $H_{t}$ is closed.

To prove closability, consider, for example, $t_{D}$, and assume
$u_{n}\underset{H}{\lra}0$ and  $(\nabla u_{n}-\nabla u_{m},\nabla 
u_{n}-\nabla
u_{m})\rightarrow0$ as $n,m\rightarrow\infty$. Then $\nabla
u_{n}\underset{H}{\lra}f$, and \be\label{e3.1} (f,\phi)=\lim_{n\to\infty}
(\nabla u_n,\phi) =-\lim_{n\to\infty}(u_n,\nabla\phi)=0, \qquad
\forall\phi\in C^\infty_0. \ee Thus $f=0$, so $t_{D}$ is closable. Similarly
one checks that $t_{N}$ and $t_{R}$ are closable. Let us denote by
$H_{2,2}^{1}(D')$ the completion of $C(\bar{D'})\bigcap
C^{\infty}(D')\bigcap H_{0}^{1}(D')$ in the norm $\Vert
u\Vert:=\left(\Vert\nabla u\Vert_{L^{2}(D')}^{2}+\Vert\nabla
u\Vert_{L^{2}(\partial D)}^{2}\right)^{1/2}$.

For an arbitrary open set $D\subset\R^{3}$ with finite volume ($|D|<\infty$,
where $|D|:=\meas D$ is the volume of $D$) the inequality \be\label{e3.2}
\|u\|_{L^3(D)} \leq \left(\|\nabla u\|_{L^2(D)} + \|u\|_{L^2(\partial
D)}\right) \ee holds, and the embedding operator
$i:H_{2,2}^{1}(D)\rightarrow L^{q}(D)$ is compact if $q<3$ and $|D|<\infty$
(\cite[p.258]{Maz}).

Consider the closed symmetric forms $t_{D}$, $t_{N}$ and $t_{R}$.
Each of these forms define a unique selfadjoint operator $A$ in $H=L^{2}(D')$,
$D(A)\subset H^{1}(D')\subset H$, $(Au,v)=t[u,v]$, $u\in D(A)$,
$v\in D[t]$, $A=A_{D}$, $A=A_{N}$, and $A=A_{R}$, respectively.

Let $L_{2,a}:=L^{2}(D',(1+|x|)^{-a})$, $a>1$, $\Vert
u\Vert_{L_{2,a}}^{2}=\int_{D'}\frac{|u|^{2}dx}{(1+|x|)^{a}}$ and $L_{0}^{2}$
be the set of $L^{2}(D')$ functions vanishing near infinity.

Fix a bounded domain $\widetilde{D}\subset D'$ whose boundary consists on
two parts: $\partial D$ and a smooth compact manifold $S$. Assume that
$i'_{1}:H^{1}(\widetilde{D}){\lra}L^{2}(\widetilde{D})$ and
$i'_{2}:H^{1}(\widetilde{D})\rightarrow L^{2}(\partial D)$ are compact. This
assumption depends only on the boundary $\partial D$ and does not depend on
the choice of the domain $\widetilde{D}\subset D'$, as long as 
$S$ is smooth.

The following theorem holds:

\begin{thm} \label{T:3.1} For any $F\in L_{0}^{2}$, each of the
boundary-value problems: \be\label{e3.3} A_iu=F, \qquad i=D, N \hbox{\ or\ }
R,\qquad A_i u=-\Delta u, \ee has a solution
$u=\lim_{0<\ep\to 0}(A-i\ep)^{-1}F:=(A-i0)^{-1}F$, $u\in
H_{loc}^{2}(D')$, $u\in L_{2,a}$, $a\in(1,2)$, and this solution is unique.  
\end{thm} Similar result holds for the operator $A-k^{2}$, where
$k=const>0$, in which case the solution $u$ satisfies the radiation
condition at infinity: \be\label{e3.4} \lim_{r\to \infty}\int_{|s|=r}|\frac
{\partial u}{\partial r}-iku|^2ds=0. \ee

\begin{proof} {[} Proof of Theorem 3.1 {]}
\textit{Uniqueness}. If $k=0$, then the uniqueness of the solution 
to (3.3) with the mentioned in Theorem 3.1 properties follows 
from the maximum
priciple. If $k>0$, then the uniqueness is established with the help
of the radiation condition as it was done in \cite{r470}, p.230. 

Since $A=A_{i}$ is selfadjoint, the equation \be\label{e3.5}
(A-i\varepsilon)u_\varepsilon=F,\qquad\varepsilon=\const>0, \qquad
A=-\Delta, \ee has a unique solution $u_{\varepsilon}\in H=L^{2}(D')$. Let
us prove that if $F\in L_{2,-a}$, then there exists the limit
\be\label{e3.6} u=\lim_{\varepsilon\to 0} u_\varepsilon, \qquad
\lim_{\varepsilon\to 0} \|u-u_\varepsilon\|_{L_{2,a}}=0, \ee and $u$ solves
\eqref{e3.3}. Thus the limiting absorption principle holds at $\lambda=0$.
Recall that the limiting absorption principle holds at a point $\lambda$ if
the limit
$$u:=\lim_{\varepsilon\rightarrow0}(A-\lambda-i\varepsilon)^{-1}F$$
exists in some sense and solves the equation $(A-\lambda)u=F$.

To prove \eqref{e3.6}, assume first that \be\label{e3.7}
\sup_{1>\varepsilon>0} \|u_\varepsilon\|_{L_{2,a}} \leq c, \ee where
$c=\const$ does not depend on $\varepsilon$. If \eqref{e3.7} holds, then
\eqref{e3.6} holds, as we will prove. Finally, we prove \eqref{e3.7}.

Let us prove that \eqref{e3.7} implies \eqref{e3.6}. Indeed,
\eqref{e3.7} implies \be\label{e3.8}
\|u_\varepsilon\|_{L^2(D'_R)}\leq c, \ee where
$D'_{R}:=D'\bigcap B_{R}$, $B_{R}:=\{ x:|x|\leq R\}$, and we
choose $R>0$ so that $\supp F\subset B_{R}$. It follows from
\eqref{e3.8} that there exists a sequence
$\varepsilon_{n}\rightarrow0$ such that
$u_{n}:=u_{\varepsilon_{n}}$ converges weakly:
$u_{n}\rightharpoonup u$ in $L^{2}(D'_{R})$. From the
relation \be\label{e3.9} t{[}u_n,\phi{]}=(F,\phi), \ee where
the form $t$ corresponds to the operator $A$ in \eqref{e3.5}
and the choice $\phi=u_{n}$ is possible, it follows that
\be\label{e3.10} t_i{[}u_n,u_n{]}=\|\nabla u_n\|_{L^2(D')}
\leq c, \qquad i=D,N \ee and \be\label{e3.11}
t_i{[}u_n,u_n{]}=\|\nabla u_n\|^2_{L^2(D')} +\int_{\partial
D} h|u_n|^2 ds\leq c, \qquad i=R. \ee From \eqref{e3.7} and
\eqref{e3.5} it follows that \be\label{e3.12} \|\Delta
u_n\|_{L^2(D'_R)} \leq c. \ee By the known elliptic
inequality: \be\label{e3.13} \|u\|_{H^2(D_1)} \leq
c(D_1,D_2) \left( \|\Delta u\|_{L^2(D_2)}+\|u\|_{L^2(D_2)}
\right), \qquad D_1\Subset D_2, \ee where $H^{2}$ is the
usual Sobolev space, it follows from \eqref{e3.10} and
\eqref{e3.8} that \be\label{e3.14} \|u_n\|_{H^2(D_1)}\leq
c,\ee where $D_{1}\Subset D'$ is any bounded strictly inner
subdomain of $D'$. By the embedding theorem, it follows that
there exists a $u$ such that \be\label{e3.15}
\lim_{n\to\infty} \|u_n-u\|_{H^1(D_1)}=0, \qquad D_1\Subset
D'. \ee 
Here and below we often use a sequence or a
subsequence denoted by the same symbol, say $u_n$, without
repeating it each time. In all cases when this is used, the
limit of any subsequence is the same and so the sequence
converges to this limit as well.

From
\eqref{e3.15} and \eqref{e3.5} it follows that
$\lim_{n\rightarrow\infty}\Vert\Delta u_{n}-\Delta u\Vert=0$, and by
\eqref{e3.13} one concludes \be\label{e3.16} \lim_{n\to\infty}
\|u_n-u\|_{H^2(D_1)} =0,\qquad D_1\Subset D'.\ee 
Passing to the limit in
\eqref{e3.5} with $\varepsilon=\varepsilon_{n}$ one gets equation
\eqref{e3.3} for $u$ in $D'$. From \eqref{e3.10} or \eqref{e3.11} it
follows that \be\label{e3.17} \hbox{\(u_n\underset{H^1(D'_R)}{\lra}u\) for
any \(R>0\).} \ee Outside the ball $B_{R}$ one has the equation
\be\label{e3.18} -\Delta u_n-i\varepsilon_n u_n=0 \hbox{\ in\ }
B'_R:=\R^3\setminus\barB_R, \qquad u_n(\infty)=0, \ee and, by Green's
formula, one gets \be\label{e3.19} u_n(x)=\int_{S_R}
\left(g_n\frac{\partial u_n}{\partial N} - u_n\frac{\partial g_n}{\partial
N}\right)ds, \quad x\in B'_R,\quad S_R:=\{s:|s|=R\}, \ee where $N$ is the
outer normal to $S_{R}$ and
$g_{n}=\frac{e^{\gamma\sqrt{\varepsilon_{n}}|x-y|}}{4\pi|x-y|}$,
$\gamma:=\frac{-1+i}{\sqrt{2}}$.

By \eqref{e3.16} and the embedding theorem, one has \be\label{e3.20}
\lim_{n\to\infty} \left( \|u_n-u\|_{L^2(S_R)} + \left\|\frac{\partial
u_n}{\partial N} -\frac{\partial u}{\partial N}\right\|_{L^2(S_R)}\right)=0.
\ee 
Passing to the limit in \eqref{e3.19} one gets \be\label{e3.21}
u(x)=\int_{S_R} \left( g\frac{\partial u}{\partial N}-u\frac{\partial
g}{\partial N}\right) ds, \qquad x\in B'_R, \qquad g:=\frac{1}{4\pi|x-y|}.
\ee Thus \be\label{e3.22} |u(x)|\leq \frac{c}{|x|}, \qquad x\in B'_R, \ee
and $u_{n}(x)$ satisfies \eqref{e3.22} with a constant $c$ independent of
$n$. 

Let $\widetilde{D}$ be a subdomain of $D'$ whose boundary has two
parts: $\partial D$ and a smooth compact manifold $S$.

If the Dirichlet condition is imposed, then the embedding 
$i':\oH1(\widetilde{D})\rightarrow L^{2}(\widetilde{D})$
is compact for any bounded domain $D$. If the Neumann condition is
imposed, then the compactness of the embedding 
$i'_{1}: H^{1}(\widetilde{D})\rightarrow L^{2}(\widetilde{D})$
imposes some restriction on the smoothness of $\partial D$, this embedding
operator is not compact for some open bounded sets $D$. This restriction
is weak: it is satisfied for any extension domain. If the Robin condition
is imposed, then we use compactness of the operators
$i'_{1}:H^{1}(\widetilde{D})\rightarrow L^{2}(\widetilde{D})$ and 
$i'_{2}:H^{1}(\widetilde{D})\rightarrow L^{2}(\widetilde{\partial D})$ 
for passing to the limit  \bee
  \lim_{n\to\infty} [(\nabla u_n,\nabla u_n)+\langle hu_n,u_n\rangle ]
=(\nabla u,\nabla u)+\langle hu,u\rangle. \eee

If the embedding operator $i'_{1}:H^{1}(\widetilde{D})\rightarrow
L^{2}(\widetilde{D})$ is compact, then \eqref{e3.10}, \eqref{e3.15} and
\eqref{e3.22} imply  
the following three conclusions: \be\label{e3.23}
\lim_{n\to\infty} \|u_n-u\|_{L^2(D'_R)}=0, \qquad \forall R<\infty, \ee
\be\label{e3.24} u_n\rightharpoonup u\hbox{\ in\ } H^1(D'_R), \ee
\be\label{e3.25} \lim_{n\to\infty} \|u_n-u\|_{L_{2,a}}=0. \ee Note that
\eqref{e3.25} follows from \eqref{e3.23} and \eqref{e3.22} if $a>1$. Indeed,
\bee
  \int_{|x|>R} \frac{|u_n-u|^2dx}{(1+|x|)^{1+a}} \leq c \int_{|x|\geq R}
\frac{dx}{(1+|x|)^{1+a}|x|^2} \leq \frac{c}{R^a}.\eee For an arbitrary small
$\delta>0$, one can choose $R$ so that $\frac{c}{R^{a}}<\delta$ and fix such
an $R$. For a fixed $R$ one takes $n$ sufficiently large and use
\eqref{e3.23} to get \bee \int_{D'_R}
\frac{|u_n-u|^2dx}{(1+|x|)^{1+a}}<\delta. \eee This implies \eqref{e3.25}.

The limit $u$ solves problem \eqref{e3.3}. We have already proved
uniqueness of its solution. therefore not only the subsequence $u_{n}$
converges to $u$, but also $u_{\varepsilon}\rightarrow u$ as 
$\varepsilon\rightarrow0$.
We have proved that \eqref{e3.7} implies \eqref{e3.6}.

To complete the proof of the existence of the solution to \eqref{e3.3} one
has to prove \eqref{e3.7}. Suppose \eqref{e3.7} is wrong. Then there is a
sequence $\varepsilon_{n}\rightarrow0$ such that $\Vert
u_{\varepsilon_{N}}\Vert_{L_{2,a}}:=\Vert
u_{n}\Vert_{L_{2,a}}\rightarrow\infty$. Let $v_{n}:=\frac{u_{n}}{\Vert
u_{n}\Vert_{2,a}}$. Then \be\label{e3.26} \|v_n\|_{2,a}=1 \ee
\be\label{e3.27} Av_n-i\varepsilon_n v_n=\frac{F}{\|u_n\|_{2,a}}. \ee By
the above argument, 
the embedding $j: H^{1}(D')\rightarrow L_{2,a}(D')$ is compact, and 
relation \eqref{e3.26} implies the
existence of $v\in L_{2,a}$ such that \be\label{e3.28} \lim_{n\to\infty}
\|v_n-v\|_{2,a}=0,\ee
and \be\label{e3.29} A v=0.\ee 
By the uniqueness result, established above,
it follows that $v=0$. Thus \eqref{e3.28} implies
$\lim_{n\rightarrow\infty}\Vert v_{n}\Vert_{2,a}=0$. This contradicts to
\eqref{e3.26}.

This contradiction proves Theorem 3.1.
\end{proof}

\begin{rem}
The above argument is valid also for solving the problem \be\label{e3.30}
A_iu-\lambda u=F, \qquad i=D,N,R, \qquad \lambda\in\R, \ee
provided that problem \eqref{e3.30} with $F=0$ has only the trivial
solution.

One may also weaken the assumption about $F$. If $F\in L_{2,-1}$, then
\eqref{e3.21} should be replaced by \be\label{e3.31} u(x)=\int_{S_R}
\left(g\frac{\partial u}{\partial N}-u\frac{\partial g}{\partial N}\right)ds
-\int_{B'_R} g(x,y)F(y)dy. \ee If $a>3$, then, using Cauchy inequality, one
gets: \be\label{e3.32} \left| \int_{B'_R} g(x,y)F dy\right|^2 \leq c
\int_{B'_R} \frac{dy}{|x-y|^2(1+|y|)^a} \int_{B'_R} |F|^2 (1+|y|)^a dy \leq
\frac{c}{|x|^2}, \ee for large $|x|$, so that \eqref{e3.22} holds if $F\in
L_{2,-a}$, $a>3$. The rest of the argument is unchanged. \end{rem}

\begin{rem} We want to emphasize that the assumptions on the smoothness of
the boundary $\partial D$ under which we have proved existence and
uniqueness of the solutions to boundary-value problems are weaker than the
usual assumptions. Namely, for the Dirichlet condition $u=0$ on $\partial D$
no assumptions, except boundedness of $D$, are used, for the Neumann
condition $u_{N}=0$ on $\partial D$ only compactness of the embedding
operator $i'_{1}:H^{1}(\widetilde{D})\rightarrow L^{2}(\widetilde{D})$ is
used, and for the Robin boundary condition $u_{N}+hu=0$ on $\partial D$,
$0\leq h\leq c$, compactness of $i'_{1}$ and of
$i'_{2}:H^{1}(\widetilde{D})\rightarrow L^{2}(\partial D)$ is used.

Our arguments can be applied for a study of the boundary-value problems
for second-order formally selfadjoint elliptic operators and for nonselfadjoint
sectorial second-order elliptic operators. In \cite{Kat} one finds
the theory of sectorial operators and the corresponding sectorial
sesquilinear forms.
\end{rem}

 \section{Quasiisometrical mappings}

The main purpose of this section is to study boundary behavior of
quasiisometrical homeomorphisms.

\subsection{Definitions and main properties.}

 Let us start with some definitions.

\begin{defn}{\em (Quasiisometrical homeomorphisms) } Let $A$ and $B$ be two
subsets of $R^{n}$. A homeomorphism $\varphi:A\rightarrow B$ is
$Q-$quasiisometrical if for any point $x\in A$ there exists such a ball
$B(x,r)$ that \begin{equation} Q^{-1}|y-z|\leq|\varphi(y)-\varphi(z)|\leq
Q|y-z|\label{3}\end{equation}
 for any $y,z\in B(x,r)\cap A$. Here the constant $Q>0$ does not depends on
the choice of $x\in A$, but the radius $r$ may depend on $x$.  \end{defn}
Obviously the inverse homeomorphism $\varphi^{-1}:A\rightarrow B$ is also
$Q-$quasiisometrical. A homeomorphism $\varphi: A\rightarrow B$ is a
quasiisometrical homeomorphism if it is a $Q$-quaiisometrical one for some
$Q$.
 Sets $A$ and $B$ are quasiisometrically equivalent if there exists a
quasiisometrical homeomorphism $\varphi:A\rightarrow B$.

\begin{defn} {\em (Lipschitz Manifolds)} A set $M\subset R^{n}$ is an
$m$-dimensional $Q$-lipschitz manifold if for any point $a\in M$ there
exists a $Q$-quasiisometrical homeomorphism $\varphi_{a}:B(0,1)\rightarrow
R^{n}$ such that $\varphi(0)=a$ and $\varphi(B(0,1)\cap R^{m})\subset M$.
Here $R^{m}:=\{ x\in R^{n}:x_{m+1}=...=x_{n}=0\}$.  \end{defn} We are
interested in compact lipschitz manifolds that are boundaries of domains 
in
$R^{n}$ and/or in $(n-1)$-dimensional lipschitz manifolds that are dense
subsets of boundaries in the sense of $(n-1)$-dimensional Hausdorff measure
$H^{n-1}$.

\begin{defn}{\em (Class $L$)}
We call a bounded domain $U\subset R^{n}$ a domain of class $L$
if:

1. There exist a bounded smooth domain $V\subset R^{n}$ and a 
quasiisometrical
homeomorphism $\varphi:V\rightarrow U$;

2. The boundary $\partial U$ of $U$ is a $(n-1)$-dimensional lipschitz
manifold. 
\end{defn}
The following proposition is well known and will be useful for a study of
domains of class $L$ and boundary behavior of quasisisometrical 
homeomorphisms.

\begin{prop} Let $A$ and $B$ be two subsets of $R^{n}$. A homeomorphism
$\varphi:A\rightarrow B$ is $Q-$quasiisometrical if and only if for any
point $a\in A$ the following inequality holds:

\[ Q^{-1}\leq\liminf_{x\rightarrow a,x\in
A}\frac{\left|\varphi(x)-\varphi(a)\right|}
{\left|x-a\right|}\leq\limsup_{x\rightarrow
a,x\in A}\frac{\left|\varphi(x)-\varphi(a)\right|}{\left|x-a\right|}\leq
Q.\]

Here the constant $Q>0$ does not depend on the choice of $a\in U$. 
\end{prop}
This proposition is a motivation for the following definition.

\begin{defn}{\em (Quasilipschitz mappings)} \label{quasilip}Let $A$ be a set
in $R^{n}$. A mapping $\varphi:A\rightarrow R^{m}$ is $Q-$quasilipshitz
if for any $a\in A$ one has: \[ \limsup_{x\rightarrow a,x\in
A}\frac{\left|\varphi(x)-\varphi(a)\right|}{\left|x-a\right|}\leq Q\]

Here the constant $Q>0$ does not depend on the choice of $a\in A$.

A mapping is quasilipschitz if it is $Q$-quasilipschitz for some
$Q$.

A homeomorphism $\varphi:A\rightarrow B$ is a quasiisometrical homeomorphism
iff $\varphi$ and $\varphi^{-1}$ are quasilipschitz.

By definition any quasilipschitz mapping is a locally lipschitz one.
A restriction of a $Q$-quasilipshitz mapping on any subset $B\subset A$
is a quasilipschitz mapping also. 
\end{defn}

\subsection{Interior metric and boundary metrics}

Suppose $A$ is a linearly connected set in $R^{n}$. An interior
metric $\mu_{A}$ on $A$ can be defined by the following way:

\begin{defn}
For any $x,y\subset A$ \[
\mu_{A}(x,y)=\inf_{\gamma_{x,y}}l(\gamma_{x,y}),\]
 where $\gamma_{x,y}:[0,1]\rightarrow A,\:\gamma_{x,y}(0)=x,\:\gamma_{x,y}(1)=y$
is a rectifiable curve and $l(\gamma_{x,y})$ is length of $\gamma_{x,y}$. 
\end{defn}
As  follows from Definition \ref{quasilip} a $Q$-quasilipschitz
mapping can change the length of a rectifiable curve by a factor 
$Q$ at most.
Hence a $Q$-quasilipschitz mapping $\varphi:A\rightarrow B$ of a
linearly connected set $A$ onto a linearly connected set $B$ is
a lipschitz mapping of the metric space $(A,\mu_{A})$ onto the
metric space $(B,\mu_{B})$. Any $Q$-quasiisometrical homeomorphism
$\varphi:A\rightarrow B$ is a bilipschitz homeomorphism of
the metric space $(A,\mu_{A})$ onto the metric space $(B,\mu_{B})$.

Because any domain $U$ of the class $L$ is quasiisometrically equivalent to a
smooth bounded domain and for a smooth bounded domain the interior
metric is equivalent to the Euclidian metric, the interior metric
$\mu_{U}$ is equivalent to the Euclidian metric for the domain $U$
also. It means that for any domain $U \in L$
 \[
K^{-1}\left|x-y\right|\leq\mu_{U}(x,y)\leq K\left|x-y\right|\]
 for any $x,y\in U$. Here a positive constant $K$ does not depend
on the choice of the points $x,y$. Therefore for any bounded domain $U\in 
L$
any quasilipshitz mapping $\varphi:U\rightarrow V$ is a lipschitz mapping 
$\varphi:(U,\mu_{U})\rightarrow R^{m}$ for the interior metric.

We will use the following definition of locally connected domain $U\in
R^{n}$ that is equivalent to the standard one.

\begin{defn} Suppose $(x_{k}\in U),(y_{k}\in U)$ are two arbitrary
convergent sequences such that
$\liminf_{k\rightarrow\infty}\,\mu_{U}(x_{k},y_{k})>0$. Call a domain
$U\in R^{n}$ locally connected if for any such sequences one has
$\lim_{k\rightarrow\infty}x_{k}\neq\lim_{k\rightarrow\infty}y_{k}$
 
\end{defn} If a boundary of a bounded domain is a topological manifold then
this domain is locally connected. Therefore,  domains of the class $L$
are locally connected domains because their boundaries are compact 
lipschitz
manifolds.

\begin{defn}
Let $A$ be a closed linearly connected subset of $R^{n}$ and $H^{n-1}(A)>0$.
Call the interior metric $\mu_{A}$ a quasieuclidean metric almost
everwhere if there exists a closed set $Q\subset A$ with $H^{n-1}(Q)=0$,
such that for any point $x\in A\setminus Q$ the following condition holds:

There exists such ball $B(x,r)$ that for any $y,z\in B(x,r)$
\[
\frac{1}{K}|y-z|\leq\mu_{A}(y,z)\leq K|y-z|, \]
 where $K=const>0$ does not depends on choice $y,z$
and $x$.
\end{defn}

By definition of lipschitz manifolds any domain of the class $L$
is quasieuclidean at any boundary point.

\begin{defn} Suppose $U$ is a domain in $R^{n}$ and 
$x_{0},y_{0}\in\partial
U$. Let us call the following 
quantity\[
\widetilde{\mu}_{\partial
U}(x_{0},y_{0}):=\lim_{\varepsilon\rightarrow0}
\inf_{\left|x-x_{0}\right|<\varepsilon,\left|y-y_{0}
\right|<\varepsilon}\mu_{U}(x,y)\]
relative interior boundary metric.
\end{defn}

Because boundary of any domain $U$ of the class $L$ is a compact 
lipschitz manifold, the relative interior boundary metric on $\partial U$ 
is
equivalent to the interior boundary metric on $\partial U$ for such
domains. This motivates the following definition:

\begin{defn}
A bounded domain $U\subset R^{n}$ has an almost quasiisometrical
boundary if $H^{n-1}(\partial U)<\infty$ and there exists a
closed set $A\in\partial U$ with $H^{n-1}(A)=0$ such that for any point
$x_{0}\in\partial U\setminus A$ the following condition holds:

There exists a ball $B(x_{0},r)$, $B(x_{0},r)\cap A=\emptyset,$ such that 
for any
$x,y\in\partial U\cap B(x_{0},r)$ one has:  \end{defn} \[ 
\frac{1}{K}\mu_{\partial
U}(x,y)\leq \widetilde{ \mu_{\partial U}}(x,y)\leq K\mu_{\partial 
U}(x,y),\]
 where $K=const>0$ does not depend on the choice of $x_{0},x$ and $y$.

We will use for the two-sided inequalities similar to the above one 
the
following short notation $\widetilde{\mu}_{\partial U}(x,y)\sim\mu_{\partial
U}(x,y)$.

If a domain $U$ has an almost quasiisometric boundary $\partial U$ and this
boundary is locally almost quasieuclidian then $\mu_{\partial
U}(x,y)\sim\left|x-y\right|$ for any $x,y\in\partial U$.

\begin{defn}
We call a bounded domain $U\subset R^{n}$ an almost quasiisometrical
domain if $H^{n-1}(\partial U)<\infty$ and there exists such a closed
set $A\in\partial U$, with $H^{n-1}(A)=0,$ that the following condition
holds:

There exists a ball $B(x_{0},r)\cap A=\emptyset$ such that for any
$x,y\in\partial U\cap B(x_{0},r)$ one has:

\[
\mu_{\partial U}(x,y)\sim\ \widetilde{mu}_{\partial U}(x,y)\sim\left|x-y\right|.\]

\end{defn}
By the extension theorem for lipschitz mappings
any quasiisometrical homeomorphism $\varphi$ of a smooth bounded
domain in $R^{n}$ onto a domain $V$ in $R^{n}$ has a lipschitz extension
$\widetilde{\psi}$ onto $R^{n}$. Denote by $\psi$ the restriction
of a lipschitz extension $\widetilde{\psi}$ on $\partial U$. By
continuity,  the extension $\psi$ is unique.

\begin{defn} Let $U$ be a smooth domain in $R^{n}$ and $V$ be a domain in
$R^{n}$ such that $H^{n-1}(\partial V)<\infty$. A quasiisometrical 
homeomorphism
$\varphi:U\rightarrow V$ has $N^{-1}$-property on the boundary if for any
$A\in\partial V$ with $H^{n-1}(A)=0$ one has  
$H^{n-1}(\psi^{-1}(A))=0$.  
\end{defn} The definition makes sense because the extension $\psi$ of a
quasiisometrical homeomorphism $\varphi$ on $\partial U$ is unique.

\begin{defn}{\em (Class $QI$)}
Let us call a bounded domain $V$ a domain of class $QI$ if:

1) There exists a quasisiometrical homeomorphism $\varphi:U\rightarrow V$
of a smooth bounded domain $U$ onto the domain $V$ that has
the $N^{-1}$-property on the boundary.

2) there exists such a closed set $A\in\partial V$, $H^{n-1}(A)=0$,
that $\partial V\setminus A$ is a $Q$-lipschitz manifold for some
$Q$;

3) $V$ is a locally connected almost quasiisometrical domain. 
\end{defn}
\begin{rem*}
The class $L$ is a subclass of the class $QI$.
\end{rem*}

\subsection{Boundary behavior of quasiisometrical homeomorphisms}

\begin{prop} \label{multipl} Suppose a $Q$-quasiisometrical homeomorphism
$\varphi:U\rightarrow V$ maps a smooth bounded domain $U$ onto a locally
connected domain $V$. Then there exists an extension $\psi$ of $\varphi$ on
$\partial U$ such that $\psi (\partial U)=\partial V$ and the mapping
$\psi|_{\partial U}$ is a lipshitz mapping of multiplicity one.  
\end{prop}
\begin{proof} Because $U$ is a smooth domain, $\varphi$ is a lipschitz
mapping. By the extension theorem for lipschitz mappings there exists a
$Q$-lipschitz extension $\widetilde{\psi}:R^{n}\rightarrow R^{n}$ of
$\varphi$. This extension is not necessarily a quasiisometrical 
homeomorphism.
By continuity of $\widetilde{\psi}$ and because $\varphi:U\rightarrow V$ is
a homeomorphism we have $\widetilde{\psi}(\partial U)=\partial V$.

Suppose $\psi:=\widetilde{\varphi}|_{\partial U}$ has multiplicity more 
than
one. Then there exist two different points $x_{0}, y_{0}\in\partial U$,
$x_{0}\neq y_{0}$ such that $\psi(x_{0})=\psi(y_{0})$. Choose two sequences:  
${x_{k}\in U}$ and ${y_{k}\in U}$ 
such that 
$\lim_{k\rightarrow\infty}x_{k}=x_{0},\lim_{k\rightarrow\infty}y_{k}=y_{0}$.
Because $U$ is a smooth bounded domain the interior metric $\mu_{U}$ is
equivalent to the Euclidean metric, i.e. there exists a positive 
constant
$Q$ such that $\mu_{U}(x_{k},y_{k})\geq Q^{-1}\left|x_{k}-y_{k}\right|\geq
Q^{-1}\left|x_{0}-y_{0}\right|>0$ for all
sufficiently large $k$. The homeomorphism
$\varphi:(U,\mu_{U})\rightarrow(V,\mu_{V})$ is a bi-lipschitz homeomorphism.  
Therefore
$\liminf_{k\rightarrow\infty}\mu_{V}(\varphi(x_{k}),\varphi(y_{k}))>0$.
Because $U$ is a locally connected domain
$\psi(x_{0})=\lim_{k\rightarrow\infty}\varphi(x_{k})
\neq\lim_{k\rightarrow\infty}\varphi(y_{k})=\psi(y_{0})$.
This contradiction proves the Proposition.  
\end{proof}
For any lipschitz $m$-dimensional compact manifold $M\subset
R^{n}$ and for any lipschitz mapping $\varphi:M\rightarrow R^{n}$ the set
$\varphi(M)$ is $H^{m}$-measurable for the $m$-dimensional Hausdorff measure
$H^{m}$ and $H^{m}(\varphi(M))<\infty$.

The next theorem, dealing with area formulas,
is a particular case of the result proved in \cite{AK} and used
for domains of the class $QI$. 

Let us start with an abstract version of this theorem.

\begin{defn} Call a metric space $X$ a $H^{k}$-rectifiable metric space if
there exists such finite or countable set of lipschitz mappings
$\alpha_{i}:A_{i}\rightarrow X$ of mesurable sets $A_{i}\subset R^{k}$ into
$X$ that $H^{k}(X\setminus\bigcup_{j}\alpha_{i}(A_{i}))=0$.

By the definition of the class $QI$ a boundary $\partial U$ of any domain
$U\in QI$ is a $H^{n-1}$-rectifiable metric space. 
\end{defn}
Our next definition represents an abstract version of Jacobian for 
$H^{k}$-rectifiable metric spaces.

\begin{defn}
Let $X$ and $Y$ be $H^{k}$-rectifiable complete metric spaces and
$F:X\rightarrow Y$ be a lipschitz mapping. Call the quantity

\[
J(x,F):=\lim_{r\rightarrow0}\frac{H^{k}(F(B(x,r))}{H^{k}(B(x,r))}\]
a formal Jacobian of $F$ at a point $x$. 
\end{defn}
\begin{thm}
\label{change} Suppose $X$ and $Y$ are $H^{k}$-rectifiable complete
metric spaces and $F:X\rightarrow Y$ is a lipschitz mapping of multiplicity
one.

Then

1. Formal Jacobian $J(x,F)$ exists $H^{k}-$almost everywhere;

2. The following area formula holds:

\[
\int_{X}J(x,F)dH^{k}=\int_{F(X)}dH^{k};\]

3. For any $u\in L^{1}(Y)$\[
\int_{X}u(F(x))J(x,F)dH^{k}=\int_{F(X)}u(y)dH^{K}.\]

\end{thm}

\begin{cor}
\label{Hausdorff} If a domain $V$ belongs to the class $QI$,
and  $\varphi:V\rightarrow U$ is a $Q$-quasiisometrical homeomorphism, 
then
$H^{n-1}(\partial U)<\infty$. 
\end{cor}
\begin{proof}
Any $Q$-quasiisometrical homeomorphism $\varphi:V\rightarrow U$
of a smooth domain $V\in R^{n}$ onto a domain $U$ of the class $QI$
has a lipschitz extension $\widetilde{\psi}:R^{n}\rightarrow R^{n}$.
By definition of the class $QI$ the domain $V$ is a locally connected
domain. Hence by Proposition \ref{multipl} the $Q$-lipschitz mapping
$\psi:=\widetilde{\psi}\setminus\partial V$ has multiplicity one
and $\psi(\partial V)=\partial U$. By Theorem \ref{change}
$H^{n-1}(\partial U)<\infty$. 
\end{proof}

\section{Quasiisometrical homeomorphisms and embedding operators.}

By Corollary \ref{Hausdorff}, $H^{n-1}(\partial V)<\infty$ for any
domain $V\in QI$. Therefore we can define Banach space $L^{2}(\partial V)$
using the Hausdorff measure $H^{n-1}$.

\begin{prop} \label{L2change} Let $U$ be a smooth domain and $V\in QI$. Any
$Q$-quasiisometrical homeomorphism $\varphi:U\rightarrow V$ that has
$N^{-1}$property on the boundary induces a bounded composition operator
$\psi^{*}:L^{2}(\partial V)\rightarrow L^{2}(\partial U)$ by the rule
$\psi^{*}(u)=u\circ\psi$.  \end{prop} 
\begin{proof} Denote by $m$ the
$(n-1)$-dimensional Lebesgue measure on $\partial U$ and by $\psi$ the
extension of $\phi$ onto $\partial U$.  By Theorem \ref{change} for any
$u\in L^{2}(\partial V)$ \[ \int_{\partial 
U}\left|u(\psi(x))\right|^{2}J(x,
\psi)dm=\int_{\partial V}\left|u(y)\right|^{2}dH^{m}.\]

 Suppose that there exists such a constant $K>0$ that $J(x,\psi)\geq K^{-1}$
for almost all $x\in\partial U$ . Denote by $A\in\partial U$, with
$H^{n-1}(A)=0$, a
set of all points for which the previous inequality does not hold.
Then
\[ \left\Vert \psi^{*}u\right\Vert _{L^{2}(\partial
U)}^{2}=\int_{\partial U\setminus
A}\left|u(\psi(x))\right|^{2}dH^{n-1}=\int_{\partial U\setminus
A}\left|u(\psi(x))\right|^{2}\frac{J(x,\psi)}{J(x,\psi)}dH^{n-1}\leq\]

\[ K\int_{\partial
V\setminus\psi^{-1}(A)}\left|u(\psi(x))\right|^{2}J(x,\psi)dm=K\left\Vert
u\right\Vert _{L^{2}(\partial V)}^{2}\]

The last equality is valid because $\varphi$ has the $N^{-1}$-property on
the boundary, i.e. $m(\psi^{-1}(A))=0$. Therefore $\psi^{*}:L^{2}(\partial
V)\rightarrow L^{2}(\partial U)$ is a bounded operator and $\left\Vert
\psi^{*}\right\Vert \leq K$.

To finish the proof we have to demonstrate that $J(x,\psi)\geq K^{-1}$.
Remember that any domain of the class $QI$ has an almost quasiisometric
boundary. 

It means that we can choose such a closed subset $A\subset \partial V$, 
with $H^{n-1}(A)=0,$ 
that the following property holds:

For any $x_{0}\in\partial V\setminus A$ there exists such a ball
$B(x_{0},r)\cap A=\emptyset$ that:

1) $B(x_{0},r)\cap A=\emptyset$ and $\widetilde{\mu}_{\partial
V}(x,y)\sim\mu_{\partial V}(x,y)\sim|x-y|$ for any $x,y\in\partial V\cap
B(x_{0},r)$,

2) $B(x_{0},r)\cap\partial V$ is a lipschitz manifold.

Let $B:=\psi^{-1}(A)$. Choose a point $z_{0}\in\partial U\setminus B$ and
such a ball $B(z_{0},R)$ that relations $\widetilde{\mu}_{\partial
U}(x,y)\sim\mu_{\partial U}(x,y)\sim|x-y|$ hold for any
$x,y\in\psi(B(z_{0},R))\cap\partial U$. This is possible because $U$ is a
smooth domain.

Because $\varphi$ is a $Q$-quasiisometric, the length $|\gamma|$
of any curve $\gamma\subset V$ satisfies the estimate:
 \[
\frac{1}{Q}|\varphi(\gamma)|\leq|\gamma|\leq Q|\varphi(\gamma)|\]
 where $|\varphi(\gamma)|$ is the length of the curve $\varphi(\gamma)\in
U$. In terms of the relative interior metric $\widetilde{\mu}_{\partial U}$
it means that \[ B_{\widetilde{\mu}_{\partial
U}}(x_{0},\frac{1}{Q}R)\subset\psi(B(z_{0},R)\cap\partial V\subset
B_{\widetilde{\mu}_{\partial U}}(x_{0},QR)\]
 where $x_{0}=\psi(z_{0})$. Without loss of generality we can suppose that
$\widetilde{\mu}_{\partial U}(x,y)\sim\mu_{\partial U}(x,y)\sim|x-y|$ for
any $x,y\in B_{\widetilde{\mu}_{\partial U}}(x_{0},QR)$. Finally we obtain
\begin{equation} B(x_{0},\frac{1}{K}R)\subset\psi(B(z_{0},R)\cap\partial
V\subset B(x_{0},KR)\label{1}\end{equation}
 for some constant $K$ that depends only on $Q$ and constants in
relations $\widetilde{\mu}_{\partial U}(x,y)\sim\mu_{\partial U}(x,y)\sim|x-y|$.

We have proved the inequality $J_{\psi}(x)\geq K^{-1}$ almost everywhere on 
$\partial U$. 
\end{proof}

\subsection{Compact embedding operators for rough domains.}

It is well known that the embedding operator $H^{1}(\Omega)\rightarrow
L^{2}(\partial\Omega)$ is compact for bounded smooth domains.

We will prove compactness of  the embedding operator for the class $QI$.
Then we extend the embedding theorem to the domains that are finite
unions of the $QI$-domains. Our proof is based on the following
result: a quasiisometrical homeomorphism $\varphi:U\rightarrow V$
induces a bounded composition operator $\varphi^{*}:H^{1}(V)\rightarrow H^{1}(U)$
by the rule $\varphi^{*}(u)=u\circ\varphi$ (see, for example \cite{GRe}
or \cite{Z}).

\begin{defn}
\emph{A domain $U$ is a domain of class $Q$ if it is a finite union
of elementary domains of class $QI$.}

Let us use the following result: 
\end{defn}
\begin{thm}
(see for example \cite{GRe} or \cite{Z}) Let $U$ and $V$ be domains
in $R^{n}$. A quasiisometrical homeomorphism $\varphi:U\rightarrow V$
induces a bounded composition operator $\varphi^{*}:H^{1}(V)\Rightarrow H^{1}(U)$
by the rule $\varphi^{*}(u)=u\circ\varphi$. 
\end{thm}
Combining this Theorem with Theorem \ref{L2change}, one gets:

\begin{thm}
\label{1quasi} If $U$ is a domain of the class $QI$, then the embedding
operator $i_{U}:H^{1}(U)\rightarrow L^{2}(\partial U)$ is compact. 
\end{thm}
\begin{proof}
By definition of the class $QI$ there exist a smooth bounded
domain $V$ and a quasiisometrical homemorphism $\varphi:V\rightarrow U$.
By Proposition \ref{L2change} $\varphi$ induces a bounded composition
operator $\psi^{*}:L^{2}(\partial U)\rightarrow L^{2}(\partial V)$
by the rule $\psi^{*}(u)=u\circ\psi$. Because the embedding operator
$i_{U}:H^{1}(U)\rightarrow L^{2}(\partial U)$ is compact and the
composition operator$(\varphi^{-1})^{*}:H^{1}(V)\rightarrow H^{1}(U)$,
induced by quasiisometrical homeomorphism $\varphi$, is bounded the
embedding operator $i_{V}:H^{1}(V)\rightarrow L^{2}(\partial V)$,
$i_{U}=(\varphi^{-1})^{*}\circ i_{V}\circ\left(\overline{\varphi}\right)^{*}$
is compact.
\end{proof}

To apply this result for domains of the class $Q$ we need the following 
lemma: 

\begin{lem} \label{sob} If $U$ and $V$ are domains of the class $QI$, then
the embedding operator $H^{1}(U\cup V)\rightarrow L^{2}(\partial(U\cup V))$
is compact.  \end{lem} \begin{proof} By previous proposition operators
 $i_{U}:H^{1}(U)\rightarrow L^{2}(\partial U)$ and
$i_{V}:H^{1}(V)\rightarrow L^{2}(\partial V)$ are compact. Choose a sequence
$\{ w_{n}\}\subset H^{1}(U\cup V)$, $\Vert w_{n}\Vert_{H^{1}(U\cup V)}\leq1$
for all $n$. Let $u_{n}:=w_{n}|_{\partial U}$ and $v_{n}:=w_{n}|_{\partial 
V}$.
Then $u_{n}\in L^{2}(\partial U)$ , $v_{n}\in L^{2}(\partial V)$, $\Vert
u_{n}\Vert_{L^{2}(\partial U)}\leq\left\Vert i_{U}\right\Vert $, $\Vert
v_{n}\Vert_{L^{2}(\partial V)}\leq\left\Vert i_{V}\right\Vert $.

Because the embedding operator $H^{1}(U)\rightarrow L^{2}(\partial U)$
is compact, we can choose a subsequence $\{ u_{n_{k}}\}$ of the sequence
$\{ u_{n}\}$ which converges in $L^{2}(\partial U)$ to a function
$u_{0}\in L^{2}(\partial U)$. Because the  embedding operator
$H^{1}(V)\rightarrow L^{2}(\partial V)$ is also compact we can choose
a subsequence $\{ v_{n_{k_{m}}}\}$ of the sequence $\{ v_{n_{k}}\}$
which converges in $L^{2}(\partial V)$ to a function $v_{0}\in L^{2}(\partial V)$.
One has: $u_{0}=v_{0}$ almost everywhere in $\partial U\cap\partial V$
and the function $w_{0}(x)$ which is defined as $w_{0}(x):=u_{0}(x)$
on $\partial U\cap\partial(U\cup V)$ and $w_{0}(x):=v_{0}(x)$ on
$\partial V\cap\partial(U\cup V)$ belongs to $L^{2}(\partial(U\cup V))$.

Hence

\[ \Vert w_{n_{k_{m}}}-w_{0}\Vert_{L^{2}(\partial(U\cup V))}\leq\Vert
u_{n_{k_{m}}}-u_{0}\Vert_{L^{2}(\partial U)}+\Vert
v_{n_{k_{m}}}-v_{0}\Vert_{L^{2}(\partial V)}.\]

Therefore $\Vert w_{n_{k_{m}}}-w_{0}\Vert_{L^{2}(\partial(U\cup
V))}\rightarrow0$ for $m\rightarrow\infty$ .  \end{proof} 
From Theorem
\ref{1quasi} and Lemma \ref{sob} the main result of this section follows
immediately:

\begin{thm}
\label{2quasi} If a domain \emph{$\Omega$} belongs to class $Q$
then the embedding operator $H^{1}(\Omega)\rightarrow L^{2}(\partial\Omega)$
is compact. 
\end{thm}
\begin{proof}
Let $U$ be an elementary domain of class $Q$. By Theorem \ref{1quasi}
the embedding operator $i_{U}:H^{1}(U)\rightarrow L^{2}(\partial U)$
is compact.

Because any domain $V$ of class $Q$ is a finite union
of domains of class $QI$ the result follows from Lemma
\ref{sob}. 
\end{proof}

\subsection{Examples}

Example 5.7 shows that a domain of the class $Q$ can have unfinite number 
of connected 
boundary components.
\begin{example}
\label{component} Take two domains:

1. Let domain $U$ is a union of rectangles
$P_{k}:=\{(x_{1},x_{2}):\left|x_{1}-2^{-k}\right|<2^{-k-2};0\leq
x_{2}<2^{-k-2}\}$, $k=1,2,...$ and the square
$S:=\left(0,1\right)\times(-1,0)$;

2. , $V:=\{(x_{1},x_{2}):0<x_{1}<1;10^{-1}x_{1}\leq x_{2}<1\}$.

In the book of V.Mazya \cite{Maz} it was proved that $U$ is a domain
of the class $L$. It is obvious that $V$ is also a domain of the
class $L$. Therefore $\Omega=U\cup V$ is a domain of class $Q$.
 By Theorem \ref{2quasi} the embedding operator $H^{1}(\Omega)\Rightarrow L^{2}(\partial\Omega)$
is compact.

The boundary $\partial\Omega$ of the plane domain $\Omega$ contains
countably many connected components that are boundaries of domains\[
S_{k}:=\{(x_{1},x_{2}):\left|x_{1}-2^{-k}\right|<2^{-k-2};10^{-1}x_{1}\leq
x_{2}<2^{-k-2}\}.\]
 The boundary of the rectangle $S_{0}:=\{(x_{1},x_{2}):0<x_{1}<1;-1\leq
x_{2}<1\}$ is also a large connected component of $\partial\Omega$.

Any neighboorhood of the point $\{0,0\}$ contains countably
many connected components of $\partial\Omega$ and therefore can not
be represented as a graph of any continuous function. 
\end{example}
Higher-dimensional examples can be constructed using the rotation
of the plane domain $\Omega$ around $x_{1}$-axis.

Next, we show  that the class $QI$ contains simply-connected
domains with non-trivial singularities.

Let us describe first a construction of a new quasiisometrical
homeomorphism using a given one. Suppose that $S_{k}(x)=kx$ is a
similarity transformation (which is called below a similarity) of
$R^{n}$ with the similarity coefficient $k>0$, $S_{k_{1}}(x)=k_{1}x$
is another similarity and $\varphi:U\rightarrow V$ is a $Q-$quasiisometrical
homeomorphism. Then a composition $\psi:=S_{k}\circ\varphi\circ S_{k_{1}}$
is a $k_{1}kQ$ -quasiisometrical homeomorphism.

This remark was used in \cite{GR} for construction of an example
of a domain with {}``spiral'' boundary which is quasiisometrically
equivalent to a cube. At {}``the spiral vertex'' the boundary of
the {}``spiral'' domain is not a graph of any lipschitz function.
 Here we will show that the {}``spiral'' domain
belongs to the class $QI$.
Let us recall the example from \cite{GR}.

\begin{example}
\label{spiral} We can start with the triangle $T:=\{(s,t):0<s<1,s<t<2s\}$
because $T$ is quasiisometrically equivalent to the unit square $Q_{2}=(0,1)\times(0,1)$.
Hence we need to construct only a quasiisometrical homeomorphism $\varphi_{0}$
from $T$ into $R^{2}$.

Let $(\rho,\theta)$ be polar coordinates in the plane. Define first a
mapping $\varphi:R_{+}^{2}\rightarrow R^{2}$ as follows:
$\varphi(s,t)=(\rho(s,t),\theta(s,t))$, $\rho(s,t)=s$ ,
$\theta(s,t)=2\pi\ln\frac{t}{s^{2}}$. Here
$R_{+}^{2}:=\{(s,t)|0<s<\infty,0<t<\infty\}$. An inverse mapping can be
calculated easily:
$\varphi^{-1}(\rho,\theta)=(s(\rho,\theta),t(\rho,\theta))$,
$s(\rho,\theta)=\rho$, $t(\rho,\theta))=\rho^{2}e^{\frac{\theta}{2\pi}}$.
Therefore $\varphi$ and $\varphi_{0}:=\varphi|T$ are diffeomorphisms.

The image of the ray $t=ks,$ $s>0,k>0$, is the logarithmic spiral
$\rho=k\exp(-\frac{\theta}{2\pi})$. Hence the image $S:=\varphi(T)=\varphi_{0}(T)$
is an {}``elementary spiral'' plane domain, because $\partial T$
is a union of two logarithmic spirals $\rho=\exp(-\frac{\theta}{2\pi})$,
$\rho=2\exp(-\frac{\theta}{2\pi})$ and the segment of the circle
$\rho=1$ .

The domain $T$ is a union of countably many subdomains
$T_{n}:=\{(s,t):e^{-(n+1)}<s<e^{-(n-1)},s<t<2s\}$, $n=1,2,...$ . On the
first domain $T_{1}$ the diffeomorphism $\varphi_{1}:=\varphi|T_{1}$ is
$Q-$quasiisometrical, because $\varphi_{1}$ is the restriction on $T_{1}$ of
a diffeomorphism $\varphi$ defined in $R_{+}^{2}$ and
$\overline{T_{1}}\subset R_{+}^{2}$. We do not calculate the number $Q$.

In \cite{GR} it was proven that any diffeomorphism 
$\varphi_{n}:=\varphi|T_{n}$
that is the composition $\varphi_{n}=S_{e^{-(n-1)}}\circ\varphi_{1}\circ
S_{e^{n-1}}$ of similarities $S_{e^{-(n-1)}}$, $S_{e^{n-1}}$ and the
$Q-$quasiisometrical diffeomorphism $\varphi_{1}$ is $Q-$quasiisometrical.
Therefore the diffeomorphism $\varphi_{0}$ is also $Q-$quasiisometrical, and
the {}``elementary spiral'' domain $U=\varphi_{0}(T)$ is quasiisometrically
equivalent to the unit square.  \end{example} By construction, the 
boundary of
the domain $U:=\varphi(T)$ is smooth at any point except the point 
$\{0\}$.
This domain is a locally connected domain. The quasiisometrical
homeomorphism $\varphi$ has $N^{-1}$ property because all the 
homeomorphisms
$\varphi_{n}$ have this property. Except the point $\{0\}$ the boundary
$\partial U$ is a $Q$-lipschitz manifold. All other properties of
$QI$-domains are subject of simple direct calculations. Therefore the domain
$T$ is a $QI$-domain.

\section{Conclusions} In this section we combine the
results about elliptic boundary problems with these about embedding
operators.

The first result is a formulation of Theorem \ref{T:1.3} for 
a large concrete class of rough domains. This result follows immediately 
from Theorem \ref{T:1.3}, Theorem 3.11 \cite{GR} and Theorem
\ref{2quasi}.

\begin{thm}
 If $D$ is a domain of the class $Q$, $F\in L_{0}^{2}(D),$ and $h\geq0$ is a
piecewise-continuous bounded function on $\partial D$, $h\not\equiv0$, then
problem \eqref{e1.9} has a solution in $H^{1}(D)$, this solution is
unique, and the problem \[ [u,\phi]+\int_{\partial
D}hu\bar{\phi}ds-\lambda(u,\phi)=(F,\phi),\quad\lambda=\const\in\R\]
 is of Fredholm type. 
\end{thm}

The next result is a formulation of Theorem \ref{T:3.1} for
a large class of rough exterior domains $D'$.

Fix a bounded domain $\widetilde{D}\subset D'$ whose boundary consists of
two parts $\partial D$ and a smooth compact manifold $S$. Assume that
$\widetilde{D}$ belongs to the class $Q$. By the definition of the class 
$Q$,
this assumption holds for any choice of $\widetilde{D}$ because for the
smooth component $S$ the conditions defining the  class $Q$ hold.

Then the following Theorem follows immediately from Theorem \ref{T:3.1},
Theorem 3.11 \cite{GR} and
Theorem \ref{2quasi}.

\begin{thm}
 For any $F\in L_{0}^{2}$, each of the boundary-value
problems: \be  A_iu=F, \qquad i=D, N \hbox{\ or\ }
R,\qquad A_i u=-\Delta u, \ee
has a solution $u=\lim_{\ep\downarrow0}(A-i\ep)^{-1}F:=(A-i0)^{-1}F$,
$u\in H_{loc}^{2}(D')$, $u\in L_{2,a}$, $a\in(1,2)$, and this solution
is unique. 
\end{thm}

\end{document}